\documentclass[11pt,twoside,titleauthor]{mwart}
\usepackage{amsmath,amsthm,amssymb,amscd}

 \usepackage[dvips=false,pdftex=false,vtex=false,
            papersize={210mm,297mm},
            width=135mm,height=210mm,heightrounded,
            marginratio=1:1,bindingoffset=0mm,
            headheight=19pt,headsep=6mm]{geometry}

 \newtheorem{thm}{Theorem}[section]
 \newtheorem{prop}[thm]{Proposition}
 \newtheorem{cor}[thm]{Corollary}
 \newtheorem{lem}[thm]{Lemma}

 \theoremstyle{definition}

 \usepackage{setspace}
 \usepackage{fancyhdr}

\begin{document}

\newcommand{\h}{{H\!X}}
\newcommand{\ZZ}{{\mathbb{Z}_+}}

\renewcommand{\thefootnote}{\fnsymbol{footnote}}

\title{Compact multipliers on spaces of analytic functions}
\author{Pawe\l{} Mleczko}
\date{{\emph{\small Faculty of Mathematics and Computer Science,\\ Adam
Mickiewicz University,\\  Umultowska 87, 61-614 Pozna{\'n},\\
Poland}}\\ {\small e-mail:  {\tt pml$@$amu.edu.pl} }}

 \maketitle

 \setstretch{1.55}

\begin{abstract}
\noindent In the paper compact multiplier operators on Banach
spaces of analytic functions on the unit disk with the range in
Banach sequence lattices are studied. If the domain space $X$ is
such that $H_\infty\hookrightarrow X\hookrightarrow H_1$,
necessary and sufficient conditions for compactness are presented.
Moreover, the calculation of the Hausdorff measure of
noncompactness for diagonal operators between Banach sequence
lattices is applied to obtaining the characterization of compact
multipliers in case the domain space $X$ satisfies
$H_\infty\hookrightarrow X\hookrightarrow H_2$.
\end{abstract}

The study of  coefficients of functions satisfying certain
properties has a~long history. Spaces of analytic functions on the
unit disk are of  special interest in this context. First papers
on the topic of examination Taylor coefficients of functions from
Hardy classes went back to Hardy and Littlewood  and the beginning
of the 20th century (see e.g., \cite{hardy}). Within this topic
the multipliers operators emerge as a~natural object of studies.
In this paper we investigate compactness properties of such
operators. We pay special attention to the case of the domain
spaces satisfying the property $H_\infty\hookrightarrow
X\hookrightarrow H_2$. Let us start with the definition of
a~multipliers.

An element $\lambda = \{\lambda_n\}\in \omega$ is said to be
a~\emph{multiplier} from a~sequence space $F$ into another one $E$
if $\lambda\,x\mathrel{\mathop:}=\{\lambda_n x_n\} \in E$ for any
$x=\{x_n\} \in F$. The set of multipliers from $F$ into $E$ is
denoted by $\mathcal{M}(F, E)$. If $F$ and $E$ are Banach sequence
lattices (see the definition below), then $\mathcal{M}(F, E)$ is
a~Banach sequence lattice equipped with the norm
\[
\|\lambda\|_{\mathcal{M}(F, E)} = \sup\big\{\|\lambda\,x\|_{E}; \,
\|x\|_{F} \leq 1\big\}.
\]
Every sequence $\lambda\in\mathcal{M}(F,E)$ induces
a~\emph{diagonal operator} $D_\lambda\colon F\to E$ given by
$D_\lambda x=\{\lambda_n x_n\}$ for $x=\{x_n\}\in F$.

In the paper we study multiplier operators defined on spaces of
analytic functions. If $G$ is a~linear subspace of $H(\mathbb{U})$
-- the space of analytic functions on the unit disk $\mathbb{U}$,
then we associate with $G$ the sequence space $\widehat{G}$ of
Taylor's coefficients of all functions in $G$, i.e.
\[
\widehat{G} =\Big \{\{\hat{f}(n)\}; \,
f=\sum_{n=0}^{\infty}\hat{f}(n)u_n \in G\Big\}.
\]
Here, as usual $u_n(z)\mathrel{\mathop:}=z^n$ for $n\in\ZZ$ and
$z\in \mathbb{U}$. With any $\lambda = \{\lambda_n\} \in
\mathcal{M}(\widehat{G}, E)$, we identify a~linear {\em multiplier
operator} $M_{\lambda}\colon G \to E$ given by $M_{\lambda}f =
\{\lambda_n \hat{f}(n)\}$ for any $f=
\sum_{n=0}^{\infty}\,\hat{f}(n)u_n \in G$.  From now on we shall
not distinguish the multiplier operator $M_{\lambda}$ from the
multiplier $\lambda$ associated with this operator.

It should be pointed out that in general the problem of
characterization of multipliers even for the special spaces $G$
and $E$ is a~difficult task. An unpublished result of Fefferman
states that $\{\lambda_n\}$ is a~multiplier from the Hardy space
$H_1$ on the disc $\mathbb{U}$ into $\ell_1$ if and only if
\[
\sup_{m\in \mathbb{Z}_{+}}\, \sum_{j=1}^{\infty}\, \Bigg(
\sum_{k=jm}^{(j+1)m-1} |\lambda_k| \Bigg)^2 <\infty.
\]
For a~survey of the results on multipliers from Hardy spaces $H_p$
to $\ell_q$ for various $p$ and $q$ we refer to \cite{Osikiewicz}
and the references included therein. It's worth mentioning that
$p$-summing multipliers were studied in \cite{Alm99} in case of
$H_p$ spaces and  within the general setting of abstract Hardy
spaces in \cite{nasza}.

The paper is organized as follows. In the preliminary section we
recall fundamental definitions and describe a~solid hull of spaces
$X$ satisfying  the condition $H_\infty\hookrightarrow
X\hookrightarrow H_2$. We use this fact to characterize
multipliers from certain spaces of analytic functions to Banach
sequence lattices. Then, in section 2, we switch to the
investigation of compactness of multiplier operators. We obtain
complete characterization in case $H_\infty\hookrightarrow
X\hookrightarrow H_2$ and deliver some necessary and sufficient
conditions in the general case. The proof of the main theorem is
based on calculation of the Hausdorff measure of noncompactness
for diagonal operators between Banach sequence lattices.

Let us point out that we haven't found in the literature results
concerning compactness  of operators of  that kind even for the
most natural settings i.e., the domain being the Hardy spaces
$H_p$ and range in sequence spaces $\ell_q$. Hence as corollaries
we obtain new theorems for the classical case. Note that in
\cite{buckley} the compactness of multipliers with the range in
$H_p$ spaces was discussed.

\section{Preliminaries}

We shall use standard notation and notions from Banach theory, as
presented e.g.,~in  \cite{lindenstrauss}. The term `operator'
stands for a~bounded and linear mapping while `$\hookrightarrow$'
denotes continuous inclusion between Banach spaces.

Let $E$ be a~Banach space of real (or complex) sequences on the
set of nonnegative integers $\ZZ$. If $E$ is an ideal space
equipped with the monotone norm i.e., $|y_n|\leq |x_n|$ for each
$n\in\ZZ$ and $\{x_n\}\in E$ implies $\{y_n\}\in E$ with
$\|\{y_n\}\|_E\leq \|\{x_n\}\|_E$, then $E$ is said to be a~Banach
\emph{sequence lattice}. Sequence space $F$ is said to be
\emph{solid} if it is a~Banach sequence lattice. A~\emph{solid
hull} $S(F)$ of a~sequence space $F$ is the smallest solid space
containing $F$. If $G$ is a~vector space of analytic functions
then by a~solid hull of $G$ we mean a~solid hull of $\widehat{G}$.
In what follows $\{e_n\}$ stands for a~standard basis in $c_0$.

A~Banach sequence lattice $E$ on $\ZZ$ is said to be
\emph{order-continuous} if every non negative non increasing
sequence in $X$ which converges to~$0$ pointwise converges to~$0$
in the norm topology of $E$. It can be easily seen that  if $E$ is
an order-continuous Banach sequence lattice than $\{x_n\}\in E$
implies $\big\| \sum_{k=n}^\infty x_ke_k\big\|_E\to0$ as
$n\to\infty$.

Since we consider quite general situation we present examples of
special spaces for which theorems -- that will be proved below --
can be in particular used. Let $X$ be a~rearangement invariant
space on $\mathbb{T}=[0,2\pi]$. Denote by $\h$ the Banach space of
all $f\in H(\mathbb{U})$ such that
\[ \|f\|_{\h}\mathrel{\mathop:}=\sup_{0\leq r<1} \|f_r\|_X<\infty,
\]
where as usual $f_r(t)\mathrel{\mathop:}=f(re^{it})$ for $t\in
\mathbb{T}$ and $0\leq r<1$. The spaces $\h$ are abstract variants
of the classical spaces that occur in the analysis. For instance
in case $X=L_p$, $1\leq p\leq \infty$, $H\!L_p$ is the Hardy space
$H_p$. We pay particular interest to the space $H_2$. It is well
known that $f\in H_2$ if and only if $\{\hat{f}(n)\}\in\ell_2$ and
$\|f\|_{H_2}=\|\{\hat{f}(n)\}\|_{\ell_2}$. For details on Hardy
spaces see \cite{duren}. We refer the reader also to \cite{nasza}
where the abstract Hardy spaces $\h$ were studied.

In what follows we investigate the conditions a~sequence $\lambda$
must satisfied in order to the multiplier $M_\lambda$ be compact.
Recall that a~bounded operator $T\colon X\to Y$ between Banach
spaces is \emph{compact} if the image of the unit ball $B_X$ is
a~relatively compact set in $Y$. It is well known that compact
operators have so called \emph{ideal property}, i.e., if $T\colon
X\to Y$ is compact, $S\colon X_0\to X$ and $R\colon Y\to Y_0$ are
bounded operators between Banach spaces, then the composition
$S\circ T\circ R\colon X_0\to Y_0$ is compact operator.

The Proposition hereunder describes the solid hull of $X$ in the
case of $H_\infty\hookrightarrow X\hookrightarrow H_2$.

\begin{prop}\label{solidity} Let $X$ be a~Banach space such that
$H_\infty\hookrightarrow X\hookrightarrow H_2$. Then
$S(X)=\ell_2$.
\end{prop}

\begin{proof}
By the result due to  Kisliakov it follows $S(H_\infty)=\ell_2$
(see \cite{kisliakov}). Since $f=\sum_{n=0}^\infty
\hat{f}(n)u_n\in H_2$ if and only if
$\big\{\hat{f}(n)\big\}\in\ell_2$, $\ell_2$ is
 a~solid sequence space and $A\hookrightarrow B$ implies $S(A)\subset S(B)$ for any vector
spaces $A,B$, then we have
\[
\ell_2=S(H_\infty) \subset S(X)\subset S(H_2)=\ell_2.
\]
\end{proof}

We use the above Proposition to describe the space
$\mathcal{M}(X,E)$ with $H_\infty\hookrightarrow X\hookrightarrow
H_2$. We take advantage of the following  result of Anderson and
Shields (see \cite{anderson}):

\noindent\emph{If\/ $E$ is a~solid space and $A$ vector sequences
space then $\mathcal{M}(A,E)=\mathcal{M}\big(S(A),E\big)$}.

\begin{cor}\label{multipliers}
If\/ $H_\infty\hookrightarrow X\hookrightarrow H_2$ then for any
Banach sequence lattice $E$
\[
\mathcal{M}(X,E)=\mathcal{M}(\ell_2,E),
\]
and each multiplier operator satisfies factorization
\begin{equation}\label{factorization}
\begin{CD}
X @>M_\lambda>> E \\
@V{i}VV @AA{D_\lambda}A\\
H_2 @>j>>\ell_2
\end{CD}\quad .
\end{equation}
\end{cor}

Let us mention that there are general results for the description
of the space of multipliers between Banach sequence lattices. We
refer the reader to \cite{turcy}, where  the multipliers between
Orlicz sequence spaces were calculated. Further, note that if $E$
is $2$-concave, then (see \cite{defant2})
\[
\mathcal{M}(\ell_2,E)= \big(( E'^{(1/2)} )'\big)^{(2)},
\]
where given any Banach sequence lattice $E$ and $1<p<\infty$, the
Banach sequence lattice $E^{(p)}$ is defined to be a~Banach
lattice of all sequences $x=\{x_n\}$ such that $\{|x_n|^{p}\}\in
E$ equipped with the norm $\|x\|_{E^{(p)}} =
\|\{|x_n|^{p}\}\|_{E}^{1/p}$ (see \cite{lindenstrauss}). A~Banach
sequence lattice $E$ is \emph{$p$-concave}, $1\leq p<\infty$ (see
\cite{lindenstrauss}) if there exists a~constant $C>0$ such that
for any finite set  $\{x_1,\dots,x_n \}$ of elements from $E$ the
following inequality holds:
\[
\bigg( \sum_{k=1}^n \| x_k \|_E^p \bigg)^{1/p}\leq C \bigg\|\bigg(
\sum_{k=1}^n |x_k|^p\bigg)^{1/p} \bigg\|_E.
\]

In consequence if $H_\infty\hookrightarrow X\hookrightarrow H_2$
and $E$ is 2-concave then
\[
\mathcal{M}(X,E)=\big(( E'^{(1/2)} )'\big)^{(2)}.
\]

In case $X=H_p$ with $2\leq p\leq \infty$ we obtain known result
for the classical case (see \cite{kim}, cf.~\cite[Theorem
A.5]{duren} where the description of $S(H_p)$, $2\leq p<\infty$
 was gained by a~different method).
\begin{cor}
Let $2\leq p\leq \infty$. Then the following statements are true:
\begin{enumerate}
\item[\textup{(i)}] If\/ $1\leq q\leq 2$ then $\lambda\in
\mathcal{M}(H_p,\ell_q)$ if and only if $\lambda\in\ell_r$, where
$\frac{1}{r}=\frac{1}{q}-\frac{1}{2}$. \item[\textup{(ii)}] If\/
$2\leq q\leq \infty$ then $\mathcal{M}(H_p,\ell_q)=\ell_\infty$.
\end{enumerate}
\end{cor}

\section{Compact multipliers}

In this section we present results concerning compactness of
multiplier operators. We start with  Proposition
\ref{compactness5} where we give a~sufficient condition for
$M_\lambda\colon X\to E$ to be compact. Recall that if
$\|T_n-T\|\to0$ and $T_n\colon X\to Y$ are compact operators
between Banach spaces then $T\colon X\to Y$ is compact.

\begin{prop}\label{compactness5}
Let\/ $X$ be a~Banach space such that $X\hookrightarrow H_1$ and
$E$ be an order-continuous Banach sequence lattice. If $\lambda\in
E$ then the operator $M_\lambda\colon X\to E$ is compact.
\end{prop}
\begin{proof}
Assume $\lambda\mathrel{\mathop:}=\{\lambda_n\}\in E$. Since
$X\hookrightarrow H_1$ there exists $C>0$ such that for every
$n\in\mathbb{Z}_+$
\[
|\hat{f}(n)|\leq C \|f\|_{X}.
\]
This implies $\|M_\lambda f \|_E\leq C\|\lambda\|_E\|f\|_{X}$ for
any $f\in X$ and thus $M_\lambda\colon X\to E$ is bounded.

For $\lambda\in E$ denote
$\lambda^n\mathrel{\mathop:}=\sum_{k=0}^n \lambda_ke_k$. Evidently
$M_{\lambda^n}$ is compact as a~finite rank operator. For any
$f=\sum_{n=0}^\infty \hat{f}(n)u_n\in X$  by the above inequality,
we have
\begin{equation*}\label{proof1}
\big\| (M_{\lambda^n}-M_\lambda)f \big\|_E=\Big\|
\sum_{k=n+1}^\infty
 \lambda_k\hat{f}(k) e_k \Big\|_E\leq C\, \Big\| \sum_{k=n+1}^\infty
\lambda_ke_k \Big\|_{E}\; \big\|f\big\|_{X},
\end{equation*}
and whence
\begin{equation*}
\big\| (M_{\lambda^n}-M_\lambda) \big\|\leq C\, \Big\|
\sum_{k=n+1}^\infty \lambda_ke_k \Big\|_{E}.
\end{equation*}
Since $E$ is order-continuous and  $\lambda\in E$ implies $\|
\sum_{k=n+1}^\infty \lambda_ke_k\|_{E}\to 0$, we conclude by the
remark preceding the Proposition  that $M_\lambda$ is compact.
\end{proof}

We cannot expect that the  condition from the above proposition is
also necessary for the compactness of multiplier operators.
Nevertheless, we show below that the compactness of $M_\lambda$
implies $\limsup |\lambda_n|\to0$.

\begin{prop}\label{limsup1}
Let $X$ be a~Banach space such that the sequence $\{u_n\}$ is
bounded in $X$ and $E\hookrightarrow \ell_\infty$ be a~Banach
sequence lattice. If\/ $M_\lambda\colon X\to E$ is compact then
\begin{equation}\label{limsup2}
\limsup_{n\to\infty} |\lambda_n|=0.
\end{equation}
\end{prop}

\begin{proof}
Suppose that $\limsup_{n\to\infty}|\lambda_n|>0$. This implies in
particular the existence  of a~subsequence $\{k_n\}$ and
a~constant $\varepsilon>0$ such that $\inf_{n\in\mathbb{Z}_+}
|\lambda_{k_n}|\geq \varepsilon>0$. Take the  sequence
$\{u_{k_n}\}$. Since $E\hookrightarrow\ell_\infty$, then
$\|e_n\|_E\geq C$ for some $C>0$ and any $n\in\ZZ$ and  we have
\[
\big\| M_\lambda u_{k_n}-M_{\lambda} u_{k_m} \big\|_E= \big\|
|\lambda_{k_n}| e_{k_n}+|\lambda_{k_m}|e_{k_m}\big\|_E\geq
C\varepsilon,
\]
for any $n\not=m$. This completes the proof.
\end{proof}

In what follows we shall see that the condition (\ref{limsup2})
implies compactness of $M_\lambda\colon X\to E$ if $H_\infty
\hookrightarrow X\hookrightarrow  H_2$ and $\ell_2\hookrightarrow
E\hookrightarrow \ell_\infty$. The proof is based on the
Proposition \ref{measureofnoncompactness}, where we calculated the
measure of noncompactness for a~diagonal operator $D_\lambda\colon
F\to E$ induced by a~sequence $\lambda$.

Recall that for a~bounded subset $A$ of a~metric space $X$ the
\emph{Hausdorff measure of noncompactness} of $A$ is given by
\[
\alpha(A)\mathrel{\mathop:}=\big\{ \varepsilon>0; A~\textrm{ has
a~finite $\varepsilon$ net in $X$}\big\}.
\]

If\, $T\colon X\to Y$ is a~bounded operator between Banach spaces
$X$ and $Y$, then  the \emph{Hausdorff measure of noncompactness}
$\beta(T)$ of an operator $T$ is a~value
$\beta(T)\mathrel{\mathop:}=\alpha(TB_X)$, where $B_X$ is a~unit
ball of $X$. It is clear that $T$ is compact if and only if
$\beta(T)=0$. Moreover, $\beta(T)\leq \|T\|$. We recall the main
properties of the measure of noncompactness. If $A,B$ are bounded
subsets of a~Banach space $X$, then
\begin{align*}
A\subset B & \Rightarrow \alpha(A)\leq \alpha(B),\\
\alpha(A+B) &\leq \alpha(A)+\alpha(B),\\
\alpha(tA) &=|t|\alpha(A)\quad\textrm{for each }t\in\mathbb{K}.
\end{align*}
More on that topic can be found e.g., in \cite{goebel}.

\begin{lem}\label{measureofnoncompactnesslemma}
Let $E$ be an order-continuous Banach lattice and let $A$ be
a~bounded subset in $E$. Then
\begin{equation}\label{lemma-equation}
\alpha(A)=\lim_{n\to\infty}\,\sup_{x\in A} \big\| (I-P_n)x
\big\|_E,
\end{equation}
where for $n\in\mathbb{Z}_+$\textup, $P_n\colon E\to E$ is given
by $ P_n x \mathrel{\mathop:}= \sum_{k=1}^n x_ke_k$ for
$x=\{x_n\}\in E$.
\end{lem}

\begin{proof}
Observe that by our hypothesis  $E$ is order-continuous it follows
easily
\[
\sup_{n\in\ZZ} \|P_n\|<\infty.
\]
Since $A\subset P_nA+(I-P_n)A$ and $P_nA$ has finite dimension it
is clear that
\begin{align*}
\alpha(A)\leq \alpha(P_nA)+\alpha\big( (I-P_n)A \big)\leq
\alpha\big( (I-P_n)A \big)\leq\sup_{x\in A} \big\| (I-P_n)x
\big\|.
\end{align*}
In consequence
\[
\alpha(A)\leq \lim_{n\to\infty}\,\sup_{x\in A} \big\| (I-P_n)x
\big\|_E.
\]

To prove the remaining inequality let $\varepsilon>0$ and choose
$Z \mathrel{\mathop:}= \{ z^1,\dots,z^k\}$ to be
a~$(\alpha(A)+\varepsilon)$ net of~$A$. Since
\[
A\subset Z+(\alpha(A)+\varepsilon)B_E,
\]
for any $x\in A$ there exist $z\in Z$ and $b\in B_E$ such that
\[
x=z+(\alpha(A)+\varepsilon)b.
\]
From this and the triangle inequality we have
\[
\sup_{x\in A}\big\| (I-P_n)x \big\|_E\leq \sup_{1\leq i\leq k}
\big\| (I-P_n)z^i \big\|_E+(\alpha(A)+\varepsilon).
\]
Since $E$ is order-continuous we get that $\| (I-P_n)z^i
\|_E=\big\|\sum_{k=n+1}^\infty z^i_ke_k\big\|_E\to0$ as
$n\to\infty$ and in consequence
\[
\lim_{n\to\infty}\,\sup_{x\in A} \big\| (I-P_n)x \big\|_E\leq
\alpha(A)+\varepsilon.
\]
Since $\varepsilon>0$ was arbitrary the proof is completed.
\end{proof}

\begin{prop}\label{measureofnoncompactness}
Let $F$ and $E$ be order-continuous Banach sequence lattices on
$\mathbb{Z}_+$ and $\ell_1\hookrightarrow F\hookrightarrow
E\hookrightarrow\ell_\infty$. Then
\[
\beta (D_\lambda\colon F\to E) =\limsup_{n\to\infty} |\lambda_n|.
\]
\end{prop}
\begin{proof}
Since $\ell_1\hookrightarrow F\hookrightarrow E$, we have
$\mathcal{M}(F,E)=\ell_\infty$. Let $\varepsilon>0$. There exists
a~subsequence $\{\lambda_{k_n}\}$ of $\{\lambda_n\}$ such that for
all $n\in\mathbb{Z}_+$
\[
| \lambda_{k_n} |>\limsup_{n\to\infty}|\lambda_n|-\varepsilon.
\]
Denote by $\lambda^{k}$ the sequence $\{\lambda_{k_n}\}$. From
Lemma \ref{measureofnoncompactnesslemma} we get
\begin{align*}
\beta(D_\lambda)=\alpha(D_\lambda B_F) \geq \alpha
(D_{\lambda^{k}} B_F)\geq \alpha\big( \{
\lambda^{k}e_n;n\in\mathbb{Z}_+ \} \big)\geq \limsup_{n\to\infty}
|\lambda_n|-\varepsilon.
\end{align*}
Thus
\[
\beta(D_\lambda)\geq \limsup_{n\to\infty} |\lambda_n|.
\]

We shall prove the opposite inequality. Take $\varepsilon >0$ and
observe that a~set
$A\mathrel{\mathop:}=\{n;|\lambda_n|>\limsup_{n\to\infty}|\lambda_n|+\varepsilon\}$
is finite. Hence while calculating the measure of noncompactness
of $D_\lambda$ we can assume without loss of generality that
$A=\emptyset$. We have
\[
\beta(D_\lambda)=\alpha\big( D_\lambda B_F \big)\leq \| D_\lambda
\|\leq \limsup_{n\to\infty} |\lambda_n|+\varepsilon.
\]
Finally we get
\[
\beta(D_\lambda)\leq \limsup_{n\to\infty} |\lambda_n|,
\]
and the result follows.
\end{proof}

\begin{thm}\label{compactness4}
Let $H_\infty\hookrightarrow X\hookrightarrow H_2$ and
$\ell_2\hookrightarrow E\hookrightarrow\ell_\infty$ with $E$ being
order-con\-tinuous. Operator $M_\lambda\colon X\to E$ is compact
if and only if
\begin{equation}\label{limsup}
\limsup_{n\to\infty} |\lambda_n|=0.
\end{equation}
\end{thm}

\begin{proof}
Suppose (\ref{limsup}) holds. Since $M_\lambda\colon X\to E$
admits the factorization (\ref{factorization}) and
$\ell_2\hookrightarrow E\hookrightarrow\ell_\infty$ from
Proposition \ref{measureofnoncompactness} and by the ideal
property of compact operators it follows that $M_\lambda$ is
compact. To obtain the converse we use Proposition \ref{limsup1}
and the proof is done.
\end{proof}

From Theorem \ref{compactness4} we obtain the following Corollary
in $H_p$ and $\ell_q$ case.
\begin{cor}
Let\/ $2\leq p\leq\infty$ and $ 2\leq q<\infty$. The operator
$M_\lambda\colon H_p\to \ell_q$ is compact if and only if
\[
\limsup_{n\to\infty} |\lambda_n|=0.
\]
\end{cor}

Further we consider the compactness of $M_\lambda\colon X\to E$
where $X\hookrightarrow H_2$ and $E$ satisfies some extra
condition. The proof is based upon the following theorem which is
the extension of Pitt's theorem (cf.~\cite{pitt}) on compact
operators on $\ell_p$ spaces (see \cite{defant}). To state the
theorem we need some additional terminology.

We say that a~Banach sequence lattice $E$ satisfies an \emph{upper
$p$-estimate}, (resp., a~\emph{lower $p$-estimate}), if there
exists a~constant $C>0$ such that for every choice of finitely
many pairwise disjoint elements $\{x_i\}_{i=1}^n$ in $E$, we have
\begin{align*}
\bigg\| \sum_{i=1}^n x_i \bigg\|_E &\leq C\bigg( \sum_{i=1}^n
\|x_i\|^p_E \bigg)^{1/p},
 \intertext{resp.,}
\bigg\| \sum_{i=1}^n x_i \bigg\|_E &\geq \frac{1}{C}\bigg(
\sum_{i=1}^n \|x_i\|^p_E \bigg)^{1/p}.
\end{align*}

\begin{thm}[{\cite[Theorem 1]{defant}}] \label{defant}
Let $E\hookrightarrow \ell_\infty$ and
$F\hookrightarrow\ell_\infty$ be quasi-Banach sequence lattices,
$E$ with an upper $t$-estimate and $F$ with a~lower $u$-estimate.
If $t>u$, then every operator $T$ from $E$ to $F$ is compact.
\end{thm}

We refer the reader to \cite{defant} where some application of the
above theorem was shown and examples of Banach sequence lattices
satisfying lower and upper estimates were presented.

\begin{prop}\label{compactness3}
Let $H_\infty\hookrightarrow X\hookrightarrow H_2$ and
$E\hookrightarrow\ell_\infty$ be a~Banach sequence lattice
satisfying a~lower $u$-estimate with $u<2$. Then the operator
$M_\lambda\colon X\to E$ is compact if and only if $\lambda\in
M(\ell_2,E)$.
\end{prop}
\begin{proof}
Since compact linear map is in particular continuous, from
Corollary \ref{multipliers} it follows that $\lambda\in
M(\ell_2,E)$.

To prove the converse assume that $M_\lambda$ is a~bounded linear
map and observe that as in (\ref{factorization}) $M_\lambda=i\circ
j\circ D_\lambda$, where $i\colon X\to H_2$ is the inclusion map,
$j$ is the isometry  and $D_\lambda\colon \ell_2\to E$ is
a~diagonal operator induced by a~sequence $\lambda$. Since $E$
satisfies a~lower $u$-estimate with $u<2$ from Theorem
\ref{defant} it follows that $D_\lambda$ is compact. By the ideal
property of compact operators the proof is complete.
\end{proof}

We finish the paper with giving the application of the above
Proposition to the case of $H_p$ and $\ell_q$ spaces.

\begin{cor}
Let\/ $M_\lambda\colon H_p\to \ell_q$ and $1\leq q< 2\leq p\leq
\infty$. Operator $M_\lambda$ is compact if and only if\/
$\lambda\in \ell_r$ with\/ $\frac{1}{r}=\frac{1}{q}-\frac{1}{2}$.
\end{cor}

\thispagestyle{empty}

\begin{thebibliography}{99}

\bibitem{Alm99}
I.~Almasri, \emph{Absolutely Summing Multipliers on $H^p$ Spaces}.
{J. Math. Anal. and Appl.} \textbf{237}~(1999), 451-463.

\bibitem{anderson}
J. M. Anderson and A. L. Shields, \emph{Coefficient Multipliers of
Bloch Functions}. Trans. Amer. Math. Soc. \textbf{224}(2) (1976),
255-265.

\bibitem{buckley}
S.M. Buckley, M.S. Ramanujan and D. Vukoti\'c, \emph{Bounded and
compact multipliers between Bergman and Hardy spaces}. Integr.
Equ. Oper. Theory \textbf{35} (1999), 1-19.

\bibitem{goebel}
J. Bana\'s and K. Goebel, \emph{Measures of Noncompactness in
Banach Spaces}. Lecture Notes in Pure and Applied Mathematics, 60,
Marcel Dekker, New York and Basel, 1980.

\bibitem{defant2}
A.~Defant,  M.~Masty\l o~and C.~Michels, \emph{Summing inclusion
maps between symmetric sequence spaces, a~survey}, Recent Progress
in Functional Analysis, K.D. Bierstedt et al. (eds.), Elsvier Sci.
2001.

\bibitem{defant}
A. Defant, A. L\'opez Molina and M.J. Rivera Ortun, \emph{On
Pitt's Theorem for Operators between Scalar and Vector-Valued
Quasi-Banach Sequence Spaces}. Monatsh. Math. \textbf{130} (2000),
7-18.

\bibitem{turcy}
P.B. Djakov and M.S. Ramanujan, \emph{Multipliers between Orlicz
Sequence Spaces}. Turk. J. Math. \textbf{24} (2000), 313-319.

\bibitem{duren} P.~Duren, \emph{Theory of $H^p$ spaces}. Academic
Press, San Diego, 1976.

\bibitem{hardy}
G.H. Hardy and J.E. Littlewood, \emph{Some properties of
fractional integrals}, II. Math. Z. \textbf{34} (1932), 403-439.

\bibitem{kim}
Y. Kim, \emph{Coefficient multipliers of $H^p$ and $G^p$ spaces}.
Math. Japonica \textbf{30}(5) (1985), 671-679.

\bibitem{kisliakov}
S. Kisliakov, \emph{Fourier coefficients of boundary values of
analytic functions in the disc and bidisc}. Trudy Mat. Inst.
Steklov. \textbf{155} (1981), 77-91 (in Russian).

\bibitem{lindenstrauss}
J.~Lindenstrauss and  L.~Tzafriri, \emph{Classical Banach Spaces
I}. Springer-Verlag, Berlin-Heidelberg-New York, 1979.

\bibitem{nasza} M.~Masty\l o~and P.~Mleczko, \emph{Absolutely summing
multipliers on abstract Hardy spaces}. Acta Math. Sinica
\textbf{25}(1) (2009), pp.~28.

\bibitem{Osikiewicz}
B.~Osikiewicz, \emph{Multipliers of Hardy spaces}. Questiones
Math. \textbf{27}~(2004), 57-73.

\bibitem{pitt}
H.R. Pitt, \emph{A note on bilinear forms}. J. London Math. Soc.
\textbf{11} (1936), 174-180.

\end{thebibliography}
\end{document}